\documentclass{amsart}

\usepackage{epsfig,euscript,mathtools,amscd}

\def\be{\begin{equation}}
\def\ee{\end{equation}}
\def\bt{\begin{theorem}}
\def\et{\end{theorem}}
\def\bl{\begin{lemma}}
\def\el{\end{lemma}}
\def\bc{\begin{corollary}}
\def\ec{\end{corollary}}

\def\pr{\noindent{\it Proof. }}

\newtheorem{theorem}{Theorem}[section]
\newtheorem{lemma}[theorem]{Lemma}
\newtheorem{corollary}[theorem]{Corollary}

\def\C{{\mathbb C}}
\def\P{{\mathbb P}} 
\def\Z{{\mathbb Z}} 
\def\R{{\mathbb R}} 
\def\l{\label}
\def\a{\mathbf}
\def\f{\EuScript}
\def\Ker{{\rm Ker\,}} 
\def\b{\boldsymbol } 
\def\v{{\varepsilon}} 
\def\hat{\widehat}

\mathtoolsset{showonlyrefs}

\begin{document}
\title{Commuting rational functions revisited}
\author{Fedor Pakovich}
\thanks{
This research was supported by the ISF, Grants No. 1432/18}

\address{Department of Mathematics\\ 
Ben Gurion University\\
P.O.B 653, Beer Sheva, 8410501, Israel}
\email{pakovich@math.bgu.ac.il}

\begin{abstract}
Let $B$ be a rational function of degree at least two that is neither a Latt\`es map nor conjugate to $z^{\pm n}$ or $\pm T_n$. 
We provide a method for describing the set  $C_B$ consisting of all rational functions commuting with $B.$ Specifically, 
we define an equivalence relation  $\underset{B}{\sim}$ on  $ C_B$
such that the quotient  $ C_B/\underset{B}{\sim}$ possesses the structure of a  finite group $G_B$, and
describe generators of $G_B$ in terms of the fundamental group of a special graph associated with $B$.
\end{abstract}

\maketitle



\section{Introduction}
In this paper we study commuting rational functions, 
that is rational solutions of the functional equation 
\be \l{comm} B\circ X=X\circ B.\ee 
More precisely, we fix a function $B\in \C(z)$ of degree at least two and study the set $ C_B$
consisting of all $X\in \C(z)$ such that \eqref{comm} holds. 

Functional equation \eqref{comm} was investigated already by Julia \cite{j} and Fatou \cite{f}.
In particular, they showed that commuting rational functions  $X$ and $B$ of degree at least two  have  
the same Julia set  $J=J(X)=J(B)$. 
Using Poincar\'e functions, 
Julia and Fatou proved that if $X$ and $B$ have no iterate in common and $J\neq \C\P^1$,
then, up to a conjugacy, $X$ and $B$ are either powers or Chebyshev polynomials.
The assumption $J\neq \C\P^1$ was removed by Ritt \cite{r}, who used a topological-algebraic method. Ritt proved that 
solutions of \eqref{comm} having no iterate in common
reduce 
either to powers, or to Chebyshev polynomials,  
or to Latt\`es maps.
 A  proof of the Ritt theorem based on modern dynamical methods was given by 
Eremenko \cite{e2}.

All the above results assume that $X$ and $B$ have no iterate in common. However,
commuting rational functions $X$ and $B$ which {\it do} have a common iterate, that is satisfy \be \l{ins}  B^{\circ l}=X^{\circ k} \ee for some $l,k\geq 1$ also exist.
The simplest examples of such functions can be obtained by setting 
 $$X= R^{\circ l_1}, \ \ \ \ B= R^{\circ l_2}, $$ where $R$ is an arbitrary rational function and $l_1,l_2\geq 1.$
More generally, denoting by  $Aut(R)$ the group of  M\"obius transformations 
commuting with $R$, we can set 
\be \l{vot} X=\mu_1 \circ R^{\circ l_1}, \ \ \ \ B=\mu_2 \circ R^{\circ l_2},\ee where
$\mu_1$ and $\mu_2$ are elements of $Aut(R)$ commuting  between themselves. 
However, it was shown already by Ritt (\cite{r}) that commuting rational functions satisfying 
\eqref{ins} are not exhausted by functions of the form \eqref{vot}. 
Although Ritt's method provides some insight on the structure of commuting rational functions $X$ and $B$ satisfying \eqref{ins}, it does not permit to describe this class of functions in an explicit way, and Ritt concluded his paper by saying:  
``we think that the example given above makes it
conceivable that no great order may reign in this class".

Functional equation \eqref{comm} is a particular case of the functional equation   
\be \l{i1}
A\circ X=X\circ B,
\ee 
where $A$ and $B$ are rational functions of degree at least two.
In case that \eqref {i1} is satisfied for some rational function $X$ of degree at least two, 
the function $B$ is called {\it semiconjugate} to the function $A$.
Semiconjugate  rational functions were investigated  in the recent  papers \cite{semi},   \cite{arn},  \cite{rec}, \cite{lattes},  \cite{fin}. In particular, it was shown in \cite{semi} that solutions 
of \eqref {i1} satisfying  
$\C(X,B)=\C(z),$  called {\it primitive}, can be described in terms of group actions on $\C\P^1$ or $\C$, implying strong restrictions on a possible form of $A$, $B$ and $X$. 
Any solution of \eqref {i1}  reduces to a primitive one
by a certain iterative process, and the  quantitative aspects of this reduction were studied  in the paper \cite{fin}. In particular, it was shown in 
 \cite{fin}  that if a rational function $B$ is not {\it special}, that is if
$B$ is neither a Latt\`es map, nor  conjugate to $z^{\pm n}$ or  $\pm T_n$, then 
   solutions of equations 
\eqref{comm} and \eqref{i1} obey some finiteness conditions.

Specifically, regarding to equation \eqref{comm},  it was shown in \cite{fin} that if $B$ is not special, 
then there 
exist {\it finitely many} rational functions $X_1, X_2,\dots, X_r$ such that $X$ commutes with $B$ if and only if 
\be \l{rep} X=X_j\circ B^{\circ k}\ee
for some  $j,$ $1\leq j \leq r,$  and  $k\geq 0.$
Moreover, the number $r$ and the degrees of $X_j,$ $1\leq j \leq r,$ can be bounded by numbers depending on $\deg B$ only. 
Notice that this result  immediately implies the Ritt theorem. Indeed, if $X$ commutes with  $B$, then 
any iterate  $X^{\circ l}$, $l\geq 1,$ does. Thus,  by the Dirichlet box principle, 
there exist distinct $l_1,$ $l_2$ such that 
$$X^{\circ l_1}=X_j\circ B^{\circ k_1}, \ \ \ \ \ \ \ X^{\circ l_2}=X_j\circ B^{\circ k_2}$$ for the same $j$ and 
some  $k_1,$ $k_2\geq 0$.  Therefore, if, say, $l_2>l_1,$ then
\be \l{yuy} X^{\circ l_2}=X^{\circ l_1}\circ  B^{\circ k_2-k_1},\ee
implying that 
\eqref{ins} holds for $l=l_2-l_1$ and $k=k_2-k_1,$ 
since $X$ and $B$ commute.

In this paper we provide a method for describing
the set $ C_B$ for non-special $B$. For such $B$ essentially all the information about $ C_B$ provided by the Ritt method reduces to the fact that any element of 
$ C_B$ has a common iterate
 with $B$. Thus, new approaches and techniques are needed, and we develop them in this paper.   Our main results are following. First, for any non-special rational function $B$ 
we define an equivalence relation  $\underset{B}{\sim}$ on the set $ C_B$
such that the quotient  $ C_B/\underset{B}{\sim}$ possesses the structure of a  finite group $G_B$.
Second, we describe generators of this group  in terms of the fundamental group of a special graph associated with $B$,  providing a method for describing $C_B$. Finally, we calculate $G_B$ for several classes of rational functions. Notice that our method 
of describing $C_B$ reduces the problem to the easier problem of
finding all functional decompositions $F=U\circ V$ for finitely many rational functions $F$.

In more details, for a non-special rational function $B$ we define an  
equivalence relation $\underset{B}{\sim}$ on the set $ C_B$, setting
$A_1\underset{B}{\sim} A_2$ if 
$$ A_1\circ B^{\circ l_1}=A_2\circ B^{\circ l_2}$$  
for some $l_1\geq 0,$ $l_2\geq 0$, and show that  
the multiplication of classes induced by the functional composition of their representatives
provides $ C_B/\underset{B}{\sim}$ with the structure of a finite group $ G_B$. 
The group structure on $ C_B/\underset{B}{\sim}$ offers a new look at the problem of describing $C_B$, and permits to characterize properties of $C_B$ in group 
theoretic terms. 
For example, the group $G_B$ is trivial 
if and only if 
any element of $C_B$ is an iterate of $B$, while $G_B$ is isomorphic to $Aut(B)$ 
if and only if 
any element of $C_B$ can be represented in the form $X=\mu \circ B^k,$ where $\mu\in Aut(B)$ and $k\geq 0.$

We describe generators of  $ G_B$ using a special finite graph $ \Gamma_B $ defined as follows.
Let $B$ be a rational function. Say that a rational function $\hat B$ is an {\it elementary transformation} of $B$
if there exist rational functions $U$ and $V$ such that  $B=V\circ U$ and $\hat B=U\circ V$. Say that rational functions 
 $B$ and $A$ are  {\it equivalent} and write $A\sim B$  if there exists 
a chain of elementary transformations between $B$ and $A$ (this equivalence relation should not be confused
with the previous one where the subscript $B$ is used).
Since for any M\"obius transformation $\mu$ the equality
$$B=(B\circ \mu^{-1})\circ \mu$$ holds, 
 the equivalence class $[B]$ of a rational function $B$ is a union of conjugacy classes. 
Moreover, by the result of \cite{rec}, the  class $[B]$ consists of {\it finitely many} conjugacy classes, unless $B$ is a flexible Latt\`es map.  
The graph $\Gamma_B $ is defined as a multigraph whose vertices
are in a one-to-one correspondence with some fixed representatives $B_i$ of conjugacy classes in $[B]$, and whose multiple edges connecting the vertices  corresponding to $B_i$ to $B_j$ 
are in a one-to-one correspondence with
solutions of the system  
$$B_i=V\circ U, \ \ \ B_j=U\circ V$$
in rational functions. 
In these terms, the main result of the paper about the group $G_B$ 
is a construction of a group epimorphism from  the fundamental group of the graph
$\Gamma_B $ to the group $G_B.$

The paper is organized as follows. In the second section we describe the set $C_B$ in terms of elementary transformations. 
In the third section 
we define the group $ G_B$.
In the fourth and the fifth sections we define the graph $\Gamma_B $ and construct a group epimorphism from $\pi_1(\Gamma_B )$ to  
$G_B$. We also show that if $A\sim B$, then
the groups $G_A$ and $ G_B$ are isomorphic. 
Notice that this implies in particular that if $A$ 
is a rational function such that the group $Aut(A)$ is non-trivial,
then  for any rational function $B\sim A$ the group $G_B$ is also  non-trivial, even although  $Aut(B)$ can be trivial. 
In the last case, functions of degree one in $C_A$ give rise to functions of higher degree in $C_B$ through the isomorphism
$G_A\cong G_B$. 

In the sixth section we calculate the group $G_B$ for certain classes of rational functions, and consider some examples. 
Specifically, we show  that for a wide class of rational functions, which we call  {\it generically decomposable},    
$G_B$ is isomorphic to $Aut(B)$. 
We also show that for a
polynomial $B$ the group $G_B$ is metacyclic.
Finally, we  discuss in details 
the example of commuting rational functions $B$ and $X$ satisfying condition \eqref{ins}
from the paper of Ritt \cite{r}.
In particular, we calculate the group $G_B$ which turns out to be a cyclic group of order 
three. We also provide a different example of this kind.

\section{The set $ C_B$ and elementary transformations} 
Let $B$ be  a rational function of degree at least two. We 
denote by $C_B$ the set of all rational functions commuting with $B.$

\bl \l{per} The set $C_B$ is closed with respect to the operation of composition, that  
is $A_1,A_2\in C_B$ implies $A_1\circ A_2\in C_B$. 
Furthermore, if 
$A\circ U\in C_B$ and $U\in C_B,$ then  $A\in C_B$.
\el
\pr Indeed,  if $A_1,A_2\in C_B$, then 
$$A_1\circ A_2\circ B= A_1\circ B \circ A_2=B \circ A_1 \circ A_2.$$
On the other hand, if  
$A \circ U\in C_B$ and  $U\in C_B$, then
$$B\circ A \circ U=A \circ U\circ B=A\circ B \circ U,$$
implying that $$B\circ A =A\circ B.\eqno{\Box}$$

We emphasize that we allow to elements of $C_B$ to have degree one, that is to be M\"obius transformations.
All M\"obius transformations commuting with $B$ obviously 
form a group denoted by $Aut(B)$ and called the {\it symmetry group} of $B.$ Since any  $\mu\in Aut(B)$ maps 
periodic points of $B$ of order $l\geq 1$  to themselves, and any M\"obius transformation is defined by 
its values at any three points, the symmetry group of any rational function is finite. In particular, 
$Aut(B)$ is one of the five well known finite rotation groups of the sphere: $A_4,$ $S_4,$ $A_5,$ $C_n,$ $D_{2n}.$ 
Notice that the property of $\mu\in Aut(B)$ to map 
periodic points of $B$ to periodic points can be used for a practical description of 
$Aut(B)$.

Let $B$ be a rational function. A rational function $\hat B$ is called an {\it elementary transformation} of $B$
if there exist rational functions $U$ and $V$ such that  $B=V\circ U$ and $\hat B=U\circ V$. We say that rational functions 
 $B$ and $A$ are  {\it equivalent} and write $A\sim B$  if there exists 
a chain of elementary transformations between $B$ and $A$.
Since for any M\"obius transformation $\mu$ the equality
$$B=(B\circ \mu^{-1})\circ \mu$$ holds, 
the equivalence class $[B]$ of a rational function $B$ is a union of conjugacy classes.
Thus, the relation $\sim$ can be considered as a 
weaker form of the classical conjugacy relation. The equivalence class $[B]$ 
contains infinitely many conjugacy classes if and only if 
$B$ is a flexible Latt\`es map (\cite{rec}). 

The following lemma is obtained by a direct calculation (see \cite{lattes}, Lemma 3.1).

\bl \l{lem1} Let \be \l{chh} L:\ B\rightarrow B_1 \rightarrow B_2  \rightarrow \dots \rightarrow B_s\ee 
be a sequence of elementary transformations, and   $U_i,$ $V_i,$ $1\leq i \leq s,$  rational functions such that 
$$B=V_1\circ U_1, \ \ \  B_i= U_i\circ V_i, \ \ \ \ \ 1\leq i\leq s,$$ 
and
\be \l{seqs}  U_{i}\circ V_{i}=V_{i+1}\circ U_{i+1},\ \ \ 1 \leq i \leq s-1.\ee
Then  the functions
\be \l{svii} U=U_s\circ U_{s-1}\circ \dots \circ U_{1}, \ \ \ \ V=V_{1}\circ \dots \circ V_{s-1}\circ V_s\ee
make the diagram 
\be \l{bhb}
\begin{CD} 
\C\P^1 @> B>>\C\P^1 \\ 
@V U  VV @VV U   V\\ 
 \C\P^1 @> B_s>> \C\P^1
 \\ 
@V {V}  VV @VV {V}   V\\ 
 \C\P^1 @> B>> \C\P^1, 
\end{CD} 
\ee
commutative 
and 
satisfy the equalities 
$$V\circ U=B^{\circ s}, \ \ \ \ \  \ U\circ V=B_s^{\circ s}.\eqno{\Box}$$ 
\el

It follows from Lemma \ref{lem1}, that any sequence of elementary transformations  
\eqref{chh} such that $B_s=B$ gives rise to a  rational function $U$ 
commuting with $B$,
and the main result of this section states that for non-special $B$ any element of $C_B$ can be obtained in this way. 

\bt \l{xru} Let $B$ be a non-special rational function of degree at least two. Then a rational function $X$ belongs to $C_B$ if and only if there exists a sequence of elementary transformation \eqref{chh} such that 
$B_s=B$ and $X= U_s\circ U_{s-1}\circ \dots \circ U_{1}.$
\et 

The proof of Theorem \ref{xru} uses the following two lemmas which are particular cases of Lemma 2.1 in \cite{semi} and 
 Theorem 2.18 in \cite{fin}, correspondingly. For the reader convenience we provide short independent proofs. We recall that a solution $A,X,B$ of  
\eqref {i1} is called {\it primitive} if $\C(X,B)=\C(z).$  We also mention that for an arbitrary solution $A,X,B$ of  
\eqref {i1} the equality 
\be \l{opa} \deg A=\deg B\ee 
holds.

\bl \l{asd} A solution $A,X,B$ of \eqref{i1} is primitive if and only if the algebraic 
curve 
\be \l{cur} A(x)-X(y)=0\ee is irreducible.
\el 
\pr By  the L\"uroth theorem, 
 there exists a rational function $W$  such that
 $\C(X,B)= \C(W)$, implying that  the equalities 
\be \l{of} X=X'\circ W,\ \ \  B=B'\circ W  \ee hold for
some rational functions $X'$ and $B'$ with $\C(X',B')=\C(z).$
Clearly, $x=X'(t),$ $ y=B'(t)$
is a generically one-to-one
parametrization of some irreducible component \be \l{eeg} C:\,F(x,y)=0\ee of \eqref{cur}. 
Furthermore, since 
 the degree of the projection of $C$ on $x$ (resp. $y$) is equal to 
$\deg X'$ (resp. $\deg B'$) the equalities  \be \l{e1} \deg_xF=\deg B', \ \ \ \deg_y F= \deg X'\ee hold. If  $ \C(X,B)=\C(z),$ then $\deg W=1$, and it follows from equalities \eqref{of}, \eqref{e1}, and \eqref{opa} 
that the curve 
$C$ coincides with curve \eqref{cur}, implying that \eqref{cur} is irreducible.
On the other hand, if $ \C(X,B)\neq \C(z),$ then $\deg W>1$,  and equalities \eqref{of}, \eqref{e1}, and \eqref{opa} imply that 
$C$ is proper component of \eqref{cur}.  \qed

\bl \l{primsol}
Let $A,X,B$ be a primitive solution of \eqref{i1}. Then for any $l\geq 1$ the solution $A^{\circ l},X,B^{\circ l}$ is also primitive. 
\el
\pr 
The proof is by induction on $l$. For $l=1$ the lemma is trivially true.
Assume that it is true for all $k\leq l$. By Lemma \ref{asd}, this implies that 
the algebraic curve 
$$C_k: \ A^{\circ k}(x)-X(y)=0$$ is irreducible for all $k\leq l$, and 
$$R_k: \  x=X(t), \ \ \ y=B^{\circ k}(t)$$ is its generically one-to-one parametrization.

Let $P_1,$ $P_2$ be arbitrary rational functions satisfying the equality 
\be \l{sle} A^{\circ (l+1)}\circ P_1
=X\circ P_2.\ee 
Since the curve  $C_l$ is irreducible and $R_l$ 
is its generically one-to-one parametrization, the equality 
 $$A^{\circ (l+1)}\circ P_1=A^{\circ l}\circ (A\circ P_1)
=X\circ P_2$$
 implies
that 
$$A\circ P_1=X\circ W,\ \  \ P_2=B^{\circ l}\circ W$$ 
for some $W\in\C(z)$. Furthermore, since the curve $C_1$ is also irreducible, it follows  from the first of these equalities that 
$$P_1=X\circ U, \ \ \ \ W=B\circ U$$ for some  $U\in \C(z)$. Thus, 
any pair of rational functions $P_1,P_2$ satisfying \eqref{sle} has the form  
$$P_1=X\circ U, \ \ \ \ P_2=B^{\circ (l+1)}\circ U$$ for some $U\in \C(z)$.
In particular, this implies that if the equalities 
\be \l{pta} X=P_1\circ W, \ \ \ B^{\circ (l+1)}=P_2\circ W\ee hold for some $P_1,P_2,W\in \C(z)$, then $\deg W=1$, since $P_1,P_2$ in \eqref{pta} satisfy \eqref{sle}.  Therefore,  $\C(X,B^{\circ (l+1)})=\C(z)$, that is  $A^{\circ (l+1)},X,B^{\circ (l+1)}$ is a primitive solution.
\qed

\noindent{\it Proof of Theorem \ref{xru}.} 
The sufficiency follows from Lemma \ref{lem1}. In the other direction, assume that $X\in C_B$. If $X$ is a M\"obius transformation, then  the sequence 
$$B=(B\circ X^{-1})\circ X\rightarrow X\circ (B\circ X^{-1})=B$$
is as required. So, assume that $\deg X\geq 2.$  

We observe first that there exist  a sequence
\eqref{chh} and a commutative diagram 
\be 
\begin{CD} 
\C\P^1 @> B>>\C\P^1 \\ 
@V U  VV @VV U   V\\ 
 \C\P^1 @> B_s>> \C\P^1
 \\ 
@V {X_0}  VV @VV {X_0}   V\\ 
 \C\P^1 @> B>> \C\P^1 
\end{CD} 
\ee
such that $U$ is defined by \eqref{svii}, the equality
$X=X_0\circ U$ holds, and the triple 
$B,X_0, B_s$ is a primitive solution of \eqref{i1}. 
Indeed, if $B,X,B$ is a primitive solution of \eqref{i1}, we can set $U=z,$ $X_0=X,$ $B_s=B.$
Otherwise, 
 $\C(X,B)= \C(W)$ for some $W$ with $\deg W>1$, and
substituting  equalities \eqref{of} in 
\eqref{i1} we see  that the diagram 
$$
\begin{CD} 
\C\P^1 @> B>>\C\P^1 \\ 
@V W  VV @VV W   V\\ 
 \C\P^1 @> W\circ B'>> \C\P^1
 \\ 
@V {X'}  VV @VV {X'}   V\\ 
 \C\P^1 @> B>> \C\P^1, 
\end{CD} 
$$
commutes. If the solution  $B, X',W\circ B'$  of \eqref{i1} is  primitive, we are done.
Otherwise, we can apply the above transformation to this solution. Since
$\deg X'<\deg X$, it is clear that after a finite number of steps we will obtain
a sequence of elementary transformations \eqref{chh} and 
functions $U,$ $X_0,$ and  $B_s$ as required.

To prove Theorem \ref{xru} we only must show that $\deg X_0=1$. Indeed, in this case
changing $U_s$ to $X_0\circ U_s$ and $B_s$ to $X_0\circ B_s\circ X_0^{-1}$, without loss of generality we may assume 
that $X_0=z$, so that $B_s=B$ and \eqref{chh} is the sequence required. 
Assume in contrary that $\deg X_0>1$. By Lemma \ref{primsol}, for any $l\geq 1$ the triple
$B^{\circ l},X_0,B_s^{\circ l}$ is a primitive solution of \eqref{i1}. On the other hand,
by the Ritt theorem, there exist $k$ and $l$ such that equality \eqref{ins} holds. 
Thus, 
$$B^{\circ l}=X^{\circ k}=X_0\circ (U\circ X^{\circ k-1}),$$ implying that the curve 
$$(U\circ X^{\circ k-1})(x)-y=0$$ is a component of the curve
\be \l{cucu} B^{\circ l}(x)-X_0(y)=0.\ee Moreover, this component is proper because $\deg X_0>1.$
Since, by Lemma \ref{asd},  this contradicts to the fact that
$B^{\circ l},X_0,B_s^{\circ l}$ is a primitive solution of \eqref{i1}, we conclude that 
$\deg X_0=1$. \qed

\section{The group $ G_B$}

Define an equivalence relation 
$\underset{B}{\sim}$ on the set $C_B$, setting  
$A_1\underset{B}{\sim} A_2$ if 
\be \l{zek} A_1\circ B^{\circ l_1}=A_2\circ B^{\circ l_2}\ee  
for some $l_1\geq 0,$ $l_2\geq 0$ (in order to distinguish this relation with the relation $\sim$ 
introduced in the previous section we use the subscript $B$). 
It is easy to see that $\underset{B}{\sim}$ is really an equivalence relation.
Indeed, $\underset{B}{\sim}$ is clearly reflexive and symmetric. Furthermore, if equalities \eqref{zek} and 
$$ A_2\circ B^{\circ n_1}=A_3\circ B^{\circ n_2}$$ hold, 
and $n_1\geq l_2,$ then 
$$A_1\circ B^{\circ (l_1+n_1-l_2)}= A_2\circ B^{\circ n_1}=A_3\circ B^{\circ n_2},$$     
implying that $A_1\underset{B}{\sim} A_3.$ 
Similarly, if $l_2\geq n_1,$ then 
$$A_3\circ B^{\circ (n_2+l_2-n_1)}= A_2\circ B^{\circ l_2}=A_1\circ B^{\circ l_1}.$$

\bl \l{eg}  Let $\a A$ be an equivalence class of $\underset{B}{\sim}$.  
For any $n\geq 1$ the class $\a A$   contains at most one rational function of degree $n$. Furthermore, if $A_0\in \a A$ is a function of minimal possible degree, then 
any $A\in \a A$ has the form $A=A_0\circ B^{\circ l},$ $l\geq 1.$  
Alternatively, the function $A_0$ can  be described 
as a unique function in $\a A$ which is not a rational function in $B.$

\el
\pr  If $\deg A_1=\deg A_2$ in \eqref{zek}, then $l_1=l_2,$ implying that $A_1=A_2.$ 
Furthermore, if
\be \l{xom} A\circ B^{\circ l_1}=A_0\circ B^{\circ l_2}\ee
and $l_1>l_2$, then 
$$A_0=A\circ B^{\circ (l_1-l_2)},$$ implying that $\deg A<\deg A_0$ 
in contradiction with the assumption. Therefore, 
$l_1\leq l_2$, and  hence 
$$ A=A_0\circ B^{\circ (l_2-l_1)}.$$ 
 
Moreover,   $A_0$ is not a rational function in $B,$ since if
	$A_0=A'\circ B$, then $A'$ commutes with $B$  by Lemma \ref{per}, implying that	$A'\underset{B}{\sim}A_0$ and 
$\deg A'<\deg A_0.$
On the other hand, if $A$  is an other function  in the class $\a A$ that is not a rational function in $B$, then  
\eqref{xom} implies that $l_1= l_2$ and  $A=A_0$.
\qed

For a rational function $B$ we denote  by $ G_B$
the set of equivalence classes of $\underset{B}{\sim}$ on $C_B$. 
We define a binary operation on the set $ G_B$ as follows. If $\a A_1$ and $\a A_2$ are equivalence classes of $\underset{B}{\sim}$\,, and  
$A_1\in  \a A_1$ and 
$A_2\in  \a A_2$ are their representatives, then  $ \a A_1\cdot \a A_2$ is defined as  the equivalence class containing 
$A_1\circ A_2$. It is easy to see that this operation is well-defined. Indeed, assume that  
$A_1\underset{B}{\sim}A_1'$ and $A_2\underset{B}{\sim}A_2'$. Then 
$$A_1\circ B^{\circ l_1}=A_1'\circ B^{\circ l_1'}$$ and 
$$A_2\circ B^{\circ l_2}=A_2'\circ B^{\circ l_2'},$$
implying that 
\be \l{xoo} A_1\circ B^{\circ l_1}\circ A_2\circ B^{\circ l_2}=A_1'\circ B^{\circ l_1'}\circ A_2'\circ B^{\circ l_2'}.\ee
Since $A_1,A_2\in C_B$, equality \eqref{xoo}  implies that  
$$A_1\circ A_2\circ B^{\circ (l_1+l_2)}=A_1'\circ A_2'\circ B^{\circ (l_1'+l_2')},$$
and hence $$A_1\circ A_2 \underset{B}{\sim}A_1'\circ A_2'.$$

\bt The set $ G_B$ equipped with the operation $\cdot$ is a finite group.  
\et 
\pr By definition, if
$A_i\in  \a A_i,$  $1\leq i \leq 3,$ then $( \a A_1\cdot  \a A_2)\cdot  \a A_3$
and $ \a A_1\cdot ( \a A_2\cdot  \a A_3)$ are classes containing the functions 
$(A_1\circ A_2)\circ A_3$ and $A_1\circ (A_2\circ A_3)$ 
correspondingly. 
On the other hand, $$(A_1\circ A_2)\circ A_3=A_1\circ (A_2\circ A_3),$$  
since  $\circ $  is an associative  operation on the set of rational functions. 
Therefore, the classes  $( \a A_1\cdot  \a A_2)\cdot  \a A_3$
and $ \a A_1\cdot ( \a A_2\cdot  \a A_3)$ coincide, and hence 
the operation $\cdot$ 
satisfies the associativity axiom. 

Clearly, the class $\a e$ containing the function $z$ and consisting of all iterates of $B$
serves as the unit element. Moreover, for any class $\a X$ there exists a class $\a X^{-1}$ such that 
\be \l{-1} \a X \cdot \a X^{-1}=\a X^{-1}\circ \a X =\a e.\ee
Indeed, by Theorem \ref{xru}, for any $X\in \a X$ 
there exists a sequence of elementary transformation \eqref{chh} such that 
$$X= U_s\circ U_{s-1}\circ \dots \circ U_{1}.$$ Further, it follows from Lemma \ref{lem1} that  
the function
$$Y= V_s\circ V_{s-1}\circ \dots \circ V_{1}$$ belongs to $C_B$, and  the functions $X$ and $Y$ 
satisfy 
\be \l{qqaa} X\circ Y=Y\circ X=B^{\circ s}.\ee 
Therefore, condition \eqref{-1} holds for  $\a X^{-1}$
defined  as the class containing the rational function $Y$.

Finally, by the result of the paper \cite{fin} cited in the introduction, there exist 
at most finitely many rational functions $A\in C_B$ which are not rational functions in $B$, implying by Lemma \ref{eg} that the group $G_B$ is finite.    
\qed

Notice that the above proof provides a method for the actual finding $\a X^{-1}.$
On the other hand, merely the existence of the inverse element follows from the Ritt theorem. 
Indeed, since  for any $X\in \a X$ 
there exist $l,k\geq 1$ such that \eqref{ins} holds,
for any class $\a X$ there exists $k$ such that 
$ \a X^k=\a e,$ implying that \eqref{-1} holds for  $\a X^{-1}= \a X^{k-1}$. 
Notice also that the Ritt theorem by itself does not imply that the group $G_B$ is finite, although implies that 
any its element has  finite order.

For $X\in C_B$ we  will denote by $\boldsymbol X$ the element of $G_B$ corresponding to the equivalence class of 
 $\underset{B}{\sim}$ containing $X$.

\bl \l{emb} The map $\mu\rightarrow \boldsymbol{\mu}$ is a group monomorphism  from the 
group 
$Aut(B)$ to the group $G_B.$  
\el
\pr Since functions from $Aut(B)$ have degree one, it follows from Lemma \ref{eg} that
$\boldsymbol{\mu_1}=\boldsymbol{\mu_2}$ if and only if  $\mu_1=\mu_2.$ Therefore, the map $\tau:\, \mu\rightarrow \boldsymbol{\mu}$ is injective, and it is 
easy to see that $\tau$ is a homomorphism of groups. \qed

We will denote the image of $Aut(B)$ in $G_B$ under the group monomorphism $\mu\rightarrow \boldsymbol{\mu}$ by $Aut_G(B)$.

\bl \l{ch1} The following conditions are equivalent.

\begin{enumerate}
\item [{1)}] Any $X\in C_B$ has the form $X=\mu\circ B^{\circ l}$ for some $\mu\in Aut(B)$ and $l\geq 0.$ 
\item  [{2)}] Any $X\in C_B$ of degree at least two is a rational function in $B$.
\item [{3)}] The group $G_B$ coincides with $Aut_G(B)$.
\end{enumerate} 
\el
\pr 
It is easy to see that 1) and 3) are equivalent, and that 1) implies 2).
Assume now that 2) holds, and let $X\in C_B$ be a function of degree 
at least two. By the assumption, 
$X=R_1\circ B$ for some $R\in \C(z).$ Moreover, since  by Lemma \ref{per} the function $R_1$ belongs to $C_B$, using 2) again we conclude that either $R_1\in Aut(B)$, or there exists $R_2\in \C(z)$ such that 
 $R_1=R_2\circ B$ and $R_2\in C_B$. It is clear that continuing this process we  will eventually obtain a
representation  $X=\mu\circ B^l$ for some $\mu\in Aut(B)$ and $l\geq 1.$ \qed

\section{The graph $\Gamma_{B}$} 
Let $B$ be  a rational function of degree at least two. Define $\Gamma_{B}$ as a multigraph 
whose vertices are
in a one-to-one correspondence with some fixed representatives of conjugacy classes in $[B]$, and whose multiple edges connecting  vertices  corresponding to representatives $B_i$ and $B_j$ 
are in a one-to-one correspondence with
solutions of the system  
\be \l{sys1} B_i=V\circ U, \ \ \ B_j=U\circ V\ee
in rational functions.
Notice that $\Gamma_{B}$ have loops. They  correspond to
solutions of \be \l{zer} B_i=U\circ V=V\circ U.\ee

\bl \l{lc} 
The graph $\Gamma_{B}$ does not depend on the choice of 
representatives of conjugacy classes in $[B]$.
\el
\pr Indeed, for any M\"obius transformations  $\alpha$ and $\beta$, to a
solution $U,V$ of system \eqref{sys1} corresponds a solution 
\be \l{fpr} U'= \beta\circ U\circ \alpha^{-1}, \ \ \ V'=  \alpha\circ V\circ \beta^{-1}\ee 
 of the system 
\be \l{sys2} \alpha\circ B_i\circ \alpha^{-1}=V'\circ U', \ \ \ \beta\circ B_j\circ \beta^{-1}=U'\circ V'.\ee
Furthermore, it is easy to see that formulas \eqref{fpr} provide a one-to-one correspon\-dence between
solutions of  \eqref{sys1} and  \eqref{sys2}.
\qed

\bt \l{gf} Let $B$ a rational function of degree at least two. Then the graph $\Gamma_{B}$ is finite, unless $B$ is
a flexible Latt\`es map. 
\et 
\pr By the main result of the paper \cite{rec}, the class $[B]$ contains infinitely many conjugacy classes if and only if $B$ is a flexible Latt\`es map. Therefore, if $B$ is not 
such a map, the graph $\Gamma_{B}$ contains only finitely many vertices.

Let us show now that the number of edges connecting two vertices is finite. 
Recall that two decompositions \be \l{bee} B=V\circ U, \ \ \ \ \ \ B=V'\circ U'\ee of a rational function $B$ into compositions of rational functions are called {\it equivalent} if
there exists a M\"obius transformation  $\mu$ such that 
\be \l{eq} V'=V \circ \mu^{-1}, \ \ \ \ U'=\mu\circ U.\ee
It is well known that equivalence classes of decompositions of $B$ are in one-to-one correspondence with imprimitivity systems of the monodromy group $Mon(B)$ of $B$. In particular, there exist at most finitely many such classes. 
Therefore, to prove the finiteness of the number of edges adjacent to the vertices corresponding to $B_i$ and $B_j$ it is enough to show that for any fixed solution $U,V$ of \eqref{sys1}
there exist only finitely many solutions $U'$, $V'$ of \eqref{sys1}
such that  decompositions \eqref{bee} are equivalent.  Since equalities \eqref{eq} combined with the equality 
$$U\circ V=U'\circ V'$$ imply the equality 
$$U\circ V=\mu \circ U\circ V\circ \mu^{-1},$$
the last statement follows from the finiteness of the group $Aut( U\circ V).$ \qed

Since in this paper we consider only non-special rational functions $B$, the corresponding graphs  $\Gamma_{B}$ are 
always finite by Theorem \ref{gf}. Notice that  the results 
of \cite{fin} imply that the number of vertices of  $\Gamma_{B}$ can be bounded by a number 
depending on $\deg B$ only (see Remark 5.2 in \cite{fin}).  Nevertheless, there exists no  absolute bound for the number of vertices of  $\Gamma_{B}$, and it is easy to construct  rational functions $B$ of degree $n$ for which 
the graph  $\Gamma_{B}$ contains $\approx\log_2n$ vertices (see \cite{semi}, p. 1241).

We always will assume that the representative of the conjugacy class of the function $B$ 
in $\Gamma_{B}$ 
is the function $B$ itself. 
Abusing notation, below we will call the functions $B_j$ simply ``vertices'' of $\Gamma_B$. 
Notice that for each vertex $B_j$ of $\Gamma_B$ there exists at least one loop starting and ending at 
$B$ which corresponds to the solution 
\be \l{ewq0} B=B\circ z=z\circ B\ee of \eqref{sys1}.
More generally, 
the solutions 
\be \l{ewq}B=(\mu^{-1} \circ B)\circ \mu=\mu\circ (\mu^{-1} \circ B), \ \  \ \ \  \ \ \  \ \mu\in Aut(B),\ee
give rise to $\vert Aut(B)\vert $ loops.

\begin{figure}[htbp]
\begin{center}
\epsfig{file=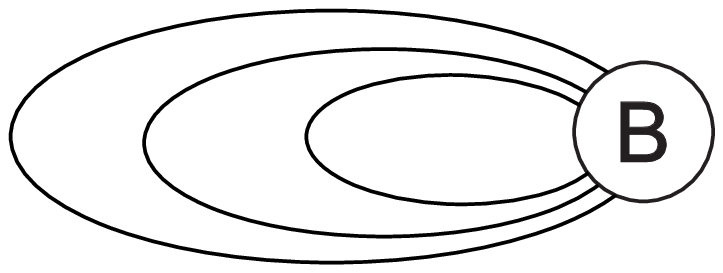,width=3cm}
\caption{ }
\label{fig0}
\end{center}
\end{figure}

\noindent{\bf Example 1.}
Assume that $B$ is an {\it indecomposable} rational function. By definition, this means that the equality $B=V\circ U$ implies that at least one of the functions $U$ and $V$ has degree one. 
In this case the equivalence class $[B]$ obviously consists of a unique conjugacy class. Thus,  
$ \Gamma_B$ has a unique vertex, and  all edges of $ \Gamma_B$ are loops corresponding 
 to solutions of
\be \l{cor} B=U\circ V=V\circ U\ee such that one of the functions $U$, $V$
 has degree one. Assuming  without loss of generality that $\deg U=1$,
we see that  
$$B\circ U=U\circ V\circ U=U\circ B,$$ implying that $U\in Aut(B).$
Therefore, 
$ \Gamma_B$ has the form shown on Fig. \ref{fig0},  
and  the number of loops of $ \Gamma_B$ is equal to $|Aut(B)|$.

\begin{figure}[htbp]
\begin{center}
\epsfig{file=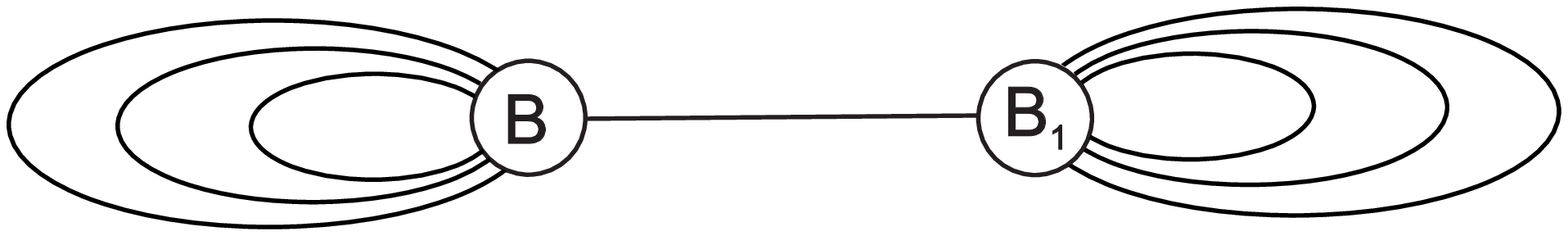,width=7.5cm}
\caption{ }
\label{fig1}
\end{center}
\end{figure}

\noindent{\bf Example 2.}
Assume now that a rational function $B$ 
has, up to equivalency \eqref{eq}, 
a  unique decomposition $B=V\circ U$ into a composition of rational functions of degree at least two, and that the same is true 
for the function $B_1=U\circ V$. 
In this case graph $  \Gamma_B$ may have  two distinct forms. 
Namely, if $B_1$ and $B$ are not conjugate, then $ \Gamma_B$ has the form shown on Fig. \ref{fig1},  
where all loops correspond to some automorphisms. 
Notice that for such $B$ and $B_1$ the groups 
$Aut(B)$  and $Aut(B_1)$ are isomorphic (see Lemma \ref{le2} below), implying that $B$ and $B_1$ have the same number of  attached loops.

On the other hand, if $B_1$ is conjugate to $B$, then without loss of generality we may assume that $B_1=B$, so that
\be \l{decc}B=V \circ U =U \circ V.\ee 
In this case the graph $ \Gamma_B$ has 
one vertex and $|Aut(B)|+1$ loops corresponding to \eqref{ewq} and \eqref{decc}. 
Notice that since by the assumption the decompositions in \eqref{decc}
are equivalent, 
the equalities 
\be \l{lis} U =V \circ \mu^{-1}, \ \ \ \ V =\mu\circ U \ee hold for some M\"obius transformation $\mu$, implying that 
$$ B=V \circ U =\mu \circ U^{\circ 2}.$$  Thus, up to a composition with a M\"obius transformation $\mu$, the function $B$ is the second iterate of some rational function $U$. Moreover, since 
$$U =V\circ\mu^{-1} =\mu\circ U \circ \mu^{-1},$$ the transfromation $\mu$ belongs to  $Aut(U).$

\noindent{\bf Example 3.}
Set $$B=-\frac{2z^2}{z^4+1}=-\frac{2}{z^2+\frac{1}{z^2}}.$$
The function $B$ is an invariant for the finite automorphism group of $\C\P^1$ generated 
by the transformations 
\be \l{alo} z\rightarrow \frac{1}{z}, \ \ \ \ z\rightarrow -z,\ee 
and its monodromy group  $Mon(B)$ is the Klein four group $\Z/2\Z\times \Z/2\Z$ having three proper imprimitivity systems.
Corresponding decompositions of 
$B$ are:
$$B=-\frac{2}{z^2-2}\circ \frac{z^2+1}{z}, \ \ \ \ \ B=-\frac{2}{z^2+2}\circ \frac{z^2-1}{z},$$
and 
\be \l{raka} B=\frac{z^2-1}{z^2+1}\circ \frac{z^2-1}{z^2+1}.\ee

\begin{figure}[htbp]
\begin{center}
\epsfig{file=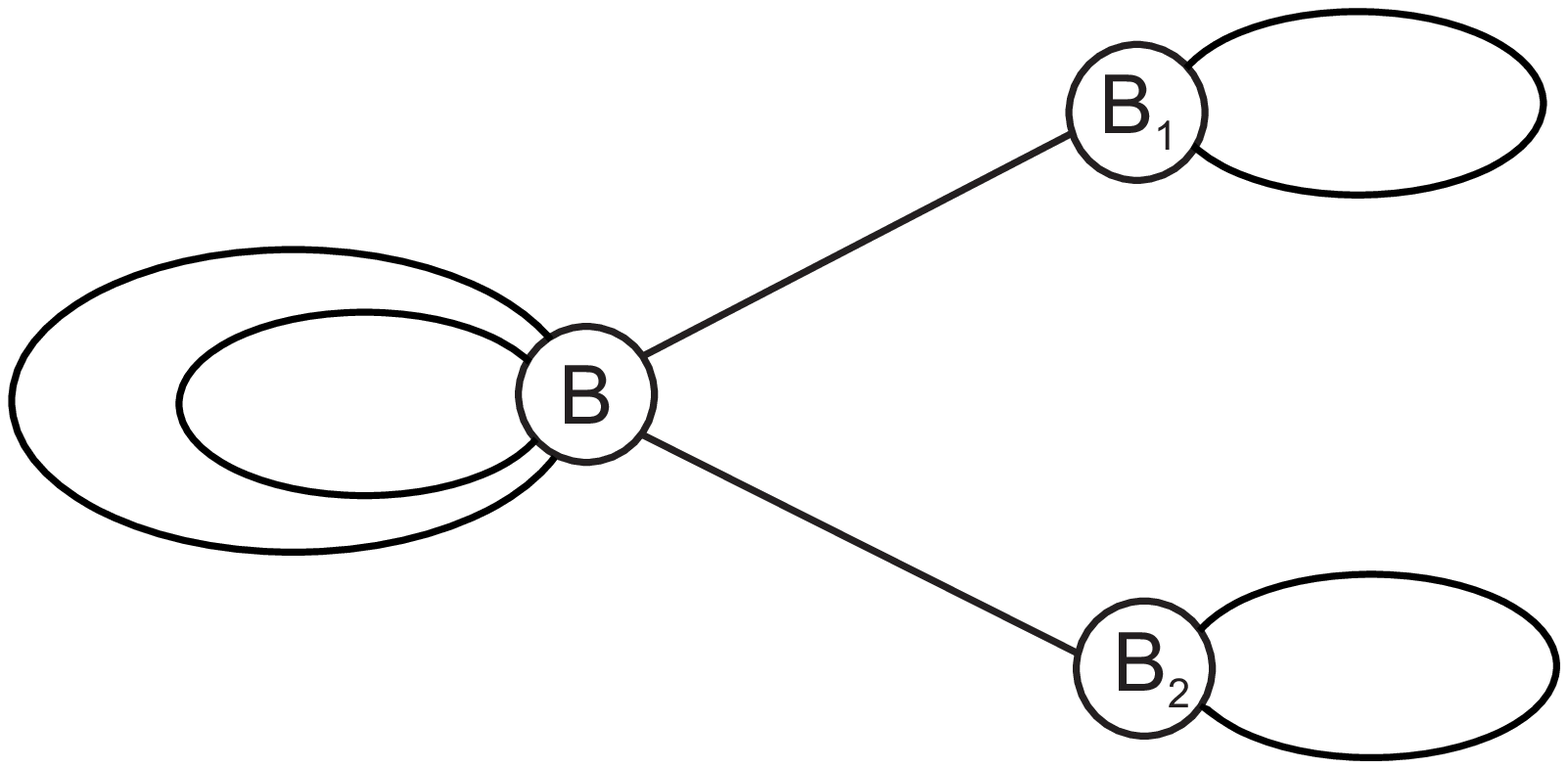,width=4.9cm}
\caption{ }
\label{fig3}
\end{center}
\end{figure}

Using for example the ``Maple'' system, one can check that 
the function 
\be \l{eg1} B_1=\frac{z^2+1}{z}\circ -\frac{2}{z^2-2}=  -\frac{1}{2}\,{\frac {{z}^{4}-4\,{z}^{2}+8}{{z}^{2}-2}}\ee 
has three critical values in $\C\P^1$, and the corresponding permutations in $Mon(B_1)$ can be identified with the  permutations $(12)(34)$, $(1243)$, and $(14)$  in $S_4.$
On the other hand, the function 
\be \l{eg2} B_2=\frac{z^2-1}{z}\circ -\frac{2}{z^2+2}=\frac{1}{2}{\frac {{z}^{2} \left( {z}^{2}+4 \right) }{{z}^{2}+2}}\ee 
has four critical values, and the corresponding permutations in $Mon(B_2)$ can be identified with 
 $(12)(34),$ $(23),$ $(12)(34)$, and $(14)$. 
Since $B_1$ and $B_2$ have a different number of critical values, they are not conjugate.
Furthermore, it is easy to see that the both groups $Mon(B_1)$ and $Mon(B_2)$  have a unique proper imprimitivity system $\{1,4\},\{2,3\}$, corresponding to decompositions \eqref{eg1} and \eqref{eg2},  implying in particular that 
$B$ is not conjugate to $B_1$ or $B_2.$ 
Finally, one can check by a direct calculations, solving the system 
$$\frac{az+b}{cz+d}\circ B=B\circ \frac{az+b}{cz+d}$$ in $a,b,c,d$, that 
the functions $B$, $B_1,$ $B_2$ have no automorphisms.
Summing up, we conclude that 
 the 
graph $ \Gamma_B$ has the form shown on Fig. \ref{fig3}.

\section{The epimorphism $\pi_1(\Gamma_B)\rightarrow G_B$} 
Considering the graph $\Gamma_B$ as a one-dimensional $CW$ complex in $\R^3$,  
we can provide each edge of $\Gamma_B$, including loops, with two opposite {\it orientations}.  
With each oriented edge $ e$ of $\Gamma_B$ we associate a rational function $\f F( e)$ as follows. 
Assume first that $e$ corresponds to solution 
\eqref{sys1} with {\it different} $B_i$ and $B_j$. Then we set $\f F( e)=U$, if the initial point of 
$ e$ is $B_i$ and the final point is $B_j$, and $\f F( e)=V$, if the orientation 
is opposite. For a loop, we simply set the value of $\f F$ equal to $U$
for one of the two corresponding oriented edges, and equal to $V$ 
 for the opposite oriented edge.
For an oriented {\it path} $$l= e_n e_{n-1}\dots  e_1$$ 
set 
$$\f F(l)=\f F( e_n)\circ \f F( e_{n-1})\circ \dots \circ \f F( e_1).$$
We emphasize  that 
since we always compose functions from right to left, we will follow this convention also for a concatenation of paths. 
Thus, a path   obtained by a concatenation of the paths $l_1$ and $l_2$ 
is denoted by 
$$l=l_2l_1,$$ and the above definition implies that 
\be \l{tov} \f F(l)=\f F(l_2)\circ \f F(l_1).\ee 
As usual, we will denote by  $l^{-1}$
the  path $l$ traversed in the opposite direction.

By construction, oriented paths from $B$ to $B_s$ correspond to sequences of elementary transformation \eqref{chh}. Furthermore, in the notation of Lemma \ref{lem1}, if $$\f F(l)=
U_s\circ U_{s-1}\circ \dots \circ U_{1},$$ then $$\f F(l^{-1})=V_{1}\circ \dots \circ V_{s-1}\circ V_s.$$
In particular, Lemma \ref{lem1} implies the following statement. 

\bl \l{sem0} 
Let $l$ be an oriented path in $\Gamma_B$ from the vertex $B$ to a vertex $B_s$ 
consisting of $k$ oriented edges.   Then
\be \l{sem} B_s\circ \f F(l)=\f F(l)\circ B ,\ee and 
\be \l{xorki} \f F(l^{-1})\circ \f F(l)= 
B^{\circ k}, \ \ \ \ \f F(l)\circ \f F(l^{-1})=
B_s^{\circ k}. \hfill\qed
\ee
\el
If $l$ is a closed path in $\Gamma_B $ starting and ending at $B$, then \eqref{sem} implies that the function $\f F(l)$ commutes with $B$, while equalities \eqref{xorki} reduce to the equalities 
\be \l{kal} \f F(l^{-1})\circ \f F(l)= \f F(l)\circ \f F(l^{-1})=
B^{\circ k}.\ee 
Thus, we obtain a map $\phi_B:\, l\rightarrow \f F(l)$
from the set of closed paths starting and ending at $B$
to the set $C_B.$

\bt \l{fg} The map $\phi_B:\, l\rightarrow \f F(l)$ descends to an epimorphism of  groups $\Phi_B:\, \pi_{1}(\Gamma_B ,B)\rightarrow G_B.$
\et 
\pr Let $\Gamma$ be a graph. Recall that an oriented path $l$ in $\Gamma$ is called {\it reduced}  if no two successive oriented edges in $l$ are opposite orientations of the same edge. Paths of the form 
$e^{-1}e$, where $e$ is an oriented edge are called {\it spurs}. 
Paths $l$ and $l'$ are called {\it equivalent} if $l'$ is obtained from $l$ by a finite number of insertions and removals of spurs between successive oriented edges or at the endpoints. In these terms, the fundamental group $\pi_{1}(\Gamma ,V)$ of the graph $\Gamma$ can be defined as the set
of equivalence classes of paths which begin and end at some fixed vertex $V$ of $\Gamma$,  
equipped with the product of classes defined in an obvious way (see e.g. Section 2.1.6 of \cite{still}).

To prove that the map $\phi_B$ descends to a map from $\pi_{1}(\Gamma_B ,B)$ to $G_B,$ we must show that whenever closed paths $l$ and $l'$ 
in $\Gamma_B $ which start and end at $B$ 
are equivalent, the rational functions $\f F(l)$ and $\f F(l')$ are in the same equivalence class of $C_B.$  Since any path is equivalent to a path with no spurs, for this purpose it is 
enough to show that if $l'$ is obtained from $l$ by an insertion  of a spur,
then $\f F(l)\underset{B}{\sim}\f F(l')$.
Assume that $$l'=l_2e^{-1}el_1,$$ where 
$l_1$ is a path from $B$ to $B_{s}$, and $l_2$ is a path from $B_{s}$ to $B$
(one of the paths $l_1$ and  $l_2$ can be empty in which case  $B_{s}=B$).
Then 
 $$\f F(l')=\f F(l_2)\circ B_{s}\circ \f F(l_1),$$ by \eqref{tov} and \eqref{kal}. It follows now from  
 \eqref{sem} that 
$$\f F(l')=\f F(l_2)\circ \f F(l_1)\circ B=\f F(l)\circ B,$$ implying that  $\f F(l)\underset{B}{\sim}\f F(l')$.
Thus, $\phi_B$ descends to a map $\Phi_B:\, \pi_{1}(\Gamma_B ,B)\rightarrow G_B,$ and \eqref{tov} implies that 
 $\Phi_B$ is a homomorphism of groups. 

Finally, it follows from Theorem \ref{xru} that $\Phi_B$ is an epimorphism. Indeed, by Theorem \ref{xru}, any $X\in C_B$ can be obtained 
from a sequence of elementary transformations \eqref{chh}. Moreover,  
we can change if necessary each of rational functions $B_i$, $1\leq i \leq s-1,$ appearing in \eqref{chh} to any desired representative of its conjugacy class, consecutively
changing the function $U_i$ to $\alpha_i\circ U_i,$ the function $B_i$ to $\alpha_i\circ B_i\circ \alpha_i^{-1},$ and the function $U_{i+1}$ to $\alpha_i^{-1}\circ U_{i+1}$ for a convenient M\"obius transformation $\alpha_i.$  
Therefore, for any $X\in C_B$ 
there exists a closed path $l$ starting and ending at $B$ such that $\f F(l)=X,$
implying that  $\Phi_B:\, \pi_{1}(\Gamma_B ,B)\rightarrow G_B$ is an epimorphism.
 \qed

\bt \l{var} Let $A$ and $B$ be equivalent rational functions. Then $G_B\cong G_A.$ 
\et
\pr Assuming that $A$ and $B$ are vertices of $\Gamma_B$, take 
 a path $s$ from $A$ to $B$ in $\Gamma_B$. Since the map 
$\psi:\,l\rightarrow s^{-1}ls$, from the set of closed paths starting and ending at $B$
to the set of closed paths starting and ending at $A$, descends to an isomorphism of 
the fundamental groups  $$\Psi:\ \pi_{1}(\Gamma_B ,B)\rightarrow \pi_{1}(\Gamma_B ,A),$$
it follows from Theorem \ref{fg} that we only must prove the equality
\be \l{rav} \Psi(\Ker \Phi_B)=\Ker \Phi_A.\ee

Let $l_0$ be a path starting and ending at $B$ such that $\f F(l_0)=B^{\circ k}$, $k\geq 1,$ 
 and let $k_0=\psi(l_0).$
Then 
$$\f F(k_0)=\f F(s^{-1})\circ \f F(l_0)\circ \f F(s)=\f F(s^{-1})\circ B^{\circ k} \circ \f F(s),$$ implying by 
\eqref{sem} and \eqref{xorki} that 
$$\f F(k_0)=\f F(s^{-1})\circ \f F(s)\circ A^{\circ k}= A^{\circ l}\circ A^{\circ k}=A^{\circ (k+l)}$$
for some $k,l\geq 1.$ This implies that 
$$\Psi(\Ker \Phi_B)\subseteq\Ker \Phi_A.$$ Similarly, considering the isomorphism inverse to $\Psi$ we obtain that 
$$\Psi^{-1}(\Ker \Phi_A)\subseteq\Ker \Phi_B.$$ This proves equality  \eqref{rav}. \qed

\section{Examples of groups $G_B$} 

\subsection{Functions with $G_B=Aut_G(B)$.}
The simplest application of Theorem \ref{fg} is the following result.

\bt \l{gr1} Let $B$ be an indecomposable non-special rational function of degree at least two. Then $G_B= Aut_G(B).$ Equivalently, $X\in C_B$ if and only if $X=\mu\circ B^l$ for some $\mu\in Aut(B)$ and $l\geq 1.$ 
\et 
\pr Since $\Gamma_B$ has a unique vertex and $|Aut(B)|$ loops corresponding to automorphisms of $B$
(see Example 1), it follows easily from Theorem \ref{fg} that 
$G_B$ is generated by $\b\mu$, $\mu\in Aut(B).$ 
Thus, $G_B= Aut_G(B).$ The second statement follows from Lemma \ref{ch1}. \qed

Notice that Theorem \ref{gr1} implies that for  a ``random'' rational function $B$ 
 the group $G_B$ is trivial, since  such a function is indecomposable and has no automorphisms. 

Theorem \ref{gr1} can be extended to a wide class of decomposable rational functions. 
Recall that a functional decomposition 
\be \l{fu} B=U_r\circ U_{r-1}\circ \dots \circ U_1\ee of a rational function $B$ is called {\it maximal}
if all $U_1,U_2,\dots, U_r$ are
indecomposable and of degree greater than one. The number $r$ is called the length of the maximal decomposition \eqref{fu}.
Two decompositions (maximal or not) having an equal number of terms
$$F=F_r\circ F_{r-1}\circ \dots \circ F_1 \ \ \ \ {\rm and} \ \ \ \  
F=G_{r}\circ G_{r-1}\circ \dots \circ G_1$$
are called equivalent if either $r=1$ and $F_1=G_1$, or $r\geq 2$ and there exist M\"obius transformations $\mu_i,$ $1\leq i \leq r-1,$ such that 
$$F_r=G_r\circ \mu_{r-1}, \ \ \ 
F_i=\mu_{i}^{-1}\circ G_i \circ \mu_{i-1}, \ \ \ 1<i< r, \ \ \ {\rm and} \ \ \ F_1=\mu_{1}^{-1}\circ G_1.
$$ Notice that all maximal decompositions of a polynomial have the same length (\cite{r1}), but this is not true 
for arbitrary rational functions (see e.g.  \cite{mp2}). 

\begin{figure}[htbp]
\begin{center}
\epsfig{file=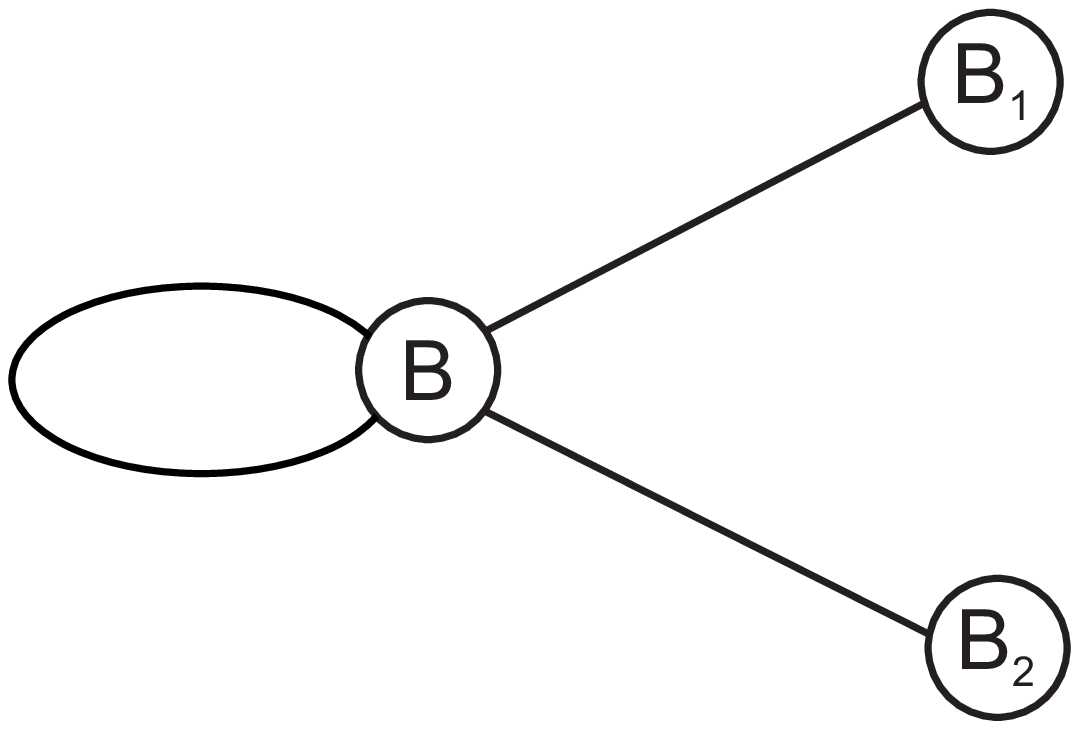,width=3.9cm}
\caption{ }
\label{fig5}
\end{center}
\end{figure}

We say that a rational function $B$ having a maximal decomposition \eqref{fu} 
is {\it generically decomposable} if the following conditions are satisfied:

\begin{itemize}

\item Each of the functions 
$$B_i= (U_{i}\circ\dots \circ U_{2}\circ U_1)\circ (U_r\circ U_{r-1}\circ \dots \circ U_{i+1}), \ \ \ \ 0\leq i \leq r-1,$$
has a unique equivalence class of maximal decompositions, 
\item The functions $B_i,$ $0\leq k \leq r-1,$ are pairwise not conjugate.
\end{itemize}
For a graph $\Gamma_B$ define  $\Gamma_B^0$ as a graph obtained from $\Gamma_B$ by removing all loops which correspond to automorphisms. For example, 
for the graph $\Gamma_B$ from Example 3 the graph $\Gamma_B^0$ is shown on Fig. \ref{fig5}.  
Recall that a {\it complete graph} is a graph in which every pair of distinct vertices is connected by a unique edge.
The complete graph on $n$ vertices is denoted by $K_n.$ 

\bl \l{le1} Assume that a non-special rational function $B$ having a maximal decomposition of length $r$
is generically decomposable. Then $\Gamma_B^0$ is the complete graph $K_r.$
\el
\pr Let \eqref{fu} be a maximal decomposition of $B.$
Since all the functions  $B_i,$ $0\leq i \leq r-1,$ are equivalent and pairwise not conjugate, the graph $\Gamma_B$ contains at least $r$ vertices.
Observe now that 
any decomposition  $B=V\circ U$ of $B$ into a composition of two rational functions of degree at least two 
has the form 
\be \l{imp} V=(U_r\circ U_{r-1}\circ \dots \circ U_{i+1})\circ \mu, \ \ \ \ \ \ U=\mu^{-1}\circ  (U_{i}\circ\dots \circ U_{2}\circ U_1),\ \ \ \ \ \ 0\leq i\leq r-1,\ee
where $\mu$ is a M\"obius transformation. Indeed, concatenating arbitrary maximal decompositions 
of $U$ and $V$ we must obtain a maximal decomposition equivalent to \eqref{fu}, implying 
that \eqref{imp} holds.
Therefore, any edge of $\Gamma_B$  adjacent to $B_0=B$ and not corresponding to an automorphism of $B$ 
is  adjacent  to one of the vertices  $B_i,$ $1\leq k \leq r-1,$ and  there exists 
exactly one edge connecting $B_0$ and $B_i$,  $1\leq k \leq r-1.$
Since the same argument holds for any $B_i,$ $0\leq k \leq r-1,$ 
we conclude that $\Gamma_B^0$ is the complete graph $K_r.$ \qed

\bl \l{le2} Assume that a non-special rational function $B$ 
is generically decomposable, and let $l$ be an oriented path from a vertex $B_{i_1}$ to a vertex $B_{i_2}$ in $\Gamma_B$. Then  for any $\mu\in Aut(B_{i_1})$ there exists $\alpha(\mu)\in Aut(B_{i_2})$
such that 
\be \l{avt} \f F(l)\circ \mu=\alpha(\mu)\circ \f F(l).\ee Furthermore, the map 
$\mu\rightarrow \alpha(\mu)$ is an isomorphism of the groups $Aut(B_{i_1})$ and  $Aut(B_{i_2}).$
 In particular, to each vertex of $\Gamma_B$ is attached the same number of loops.
\el
\pr In view of formula \eqref{tov} it is enough to prove the lemma for the case where $l$ is an oriented edge. If 
 $l$ is  a loop, then by Lemma \ref{le1} it corresponds to a solution of \eqref{sys1} of the form
\be \l{form} B_{i_1}=(\mu_0^{-1} \circ B_{i_1})\circ \mu_0=\mu_0\circ (\mu_0^{-1} \circ B_{i_1}), \ \  \ \ \  \ \ \  \ \mu_0\in Aut(B_{i_1}).\ee  Thus, either $\f F(l)=\mu_0$ or 
$\f F(l)=\mu^{-1}_0 \circ B_{i_1}$, and it is easy to see that in these cases equality \eqref{avt}  holds  for the  
automorphisms 
$$ \alpha(\mu)= \mu_0\circ \mu\circ \mu_0^{-1}, \ \ \ \  \  \alpha(\mu)= \mu_0^{-1}\circ \mu\circ \mu_0,$$
correspondingly.

 Assume now that $l$ is an oriented edge from a vertex
$B_{i_1}=V\circ U$ to a different vertex $B_{i_2}=U\circ V.$ 
Let us observe that for any $\mu\in Aut(B_{i_1})$ the decompositions $B_{i_1}=V\circ U$ and
\be \l{deco} B_{i_1}=(\mu^{-1}\circ V)\circ (U\circ \mu)\ee are equivalent, since 
 for arbitrary maximal decompositions 
of $U$ and $V$ the corresponding induced maximal decompositions of $B_{i_1}$  are equivalent.
Therefore, for any  $\mu\in Aut(B_{i_1})$
there exists a M\"obius transformation $\alpha=\alpha(\mu)$  such that 
$$\mu^{-1} \circ V=V\circ \alpha(\mu)^{-1}, \ \ \ \ U\circ \mu=\alpha(\mu)\circ U.$$
Furthermore, since 
$$B_{i_2}=U\circ V=U\circ \mu\circ \mu^{-1} \circ V=\alpha(\mu)\circ U\circ V\circ \alpha(\mu)^{-1},$$ 
the transformation  $\alpha(\mu)$ belongs to $\mu\in Aut(B_{i_2})$, 
and it is easy to see that $\mu\rightarrow \alpha(\mu)$ is a group homomorphism from $Aut(B_{i_1})$ to  $Aut(B_{i_2}).$

Finally, if $$ \nu\rightarrow \beta(\nu)$$ is a homomorphism from $Aut(B_{i_2})$ to  $Aut(B_{i_1}),$
defined by the conditions 
$$\nu^{-1} \circ U=U\circ \beta(\nu)^{-1}, \ \ \ \ V\circ \nu=\beta(\nu)\circ V,$$ and $\mu\in Aut(B_{i_1})$,
then 
$$V\circ U\circ \mu=V\circ \alpha(\mu) \circ U=\beta(\alpha(\mu)) \circ V\circ U.$$ 
Since $$V\circ U\circ \mu=\mu\circ V\circ U,$$ this implies that
$\beta\circ\alpha$ is the identical mapping of  $Aut(B_{i_1})$, and hence $\mu\rightarrow \alpha(\mu)$ is an isomorphism. \qed

\bt \l{tde} Let $B$  be a non-special generically decomposable rational function. Then  $G_B=Aut_G(B)$. 
Equivalently, $X\in C_B$ if and only if $X=\mu\circ B^l$ for some $\mu\in Aut(B)$ and $l\geq 1.$
\et 
\pr Let \eqref{fu} be a maximal decomposition of $B$. For convenience, 
define rational functions $U_i$ for $i\geq r$ setting $U_i=U_{i'}$, where $i\equiv i'\, \mod r.$
Let us recall that any decomposition $B=V\circ U$, where $U$ and $V$ are functions of degree at least two, has the form \eqref{imp}, and a similar statement holds for 
all $B_i,$ $0\leq i \leq r-1.$ Therefore, for the oriented edge $ e$ from a vertex $B_{i_1}$ to a different vertex $B_{i_2}$ the equality 
$$\f F( e)=U_{i_2}\circ\dots \circ U_{i_1+2}\circ U_{i_1+1}$$ holds, implying inductively by \eqref{tov} that 
for an arbitrary path  $l$  with no loops  from $B_{i_1}$ to $B_{i_2}$ the equality 
$$\f F(l)=U_{i_2+rk}\circ\dots \circ U_{i_1+2}\circ U_{i_1+1}=B_{i_2}^{\circ k}\circ U_{i_2}\circ\dots \circ U_{i_1+2}\circ U_{i_1+1}$$ holds
for some $k\geq 1.$ In particular, if $l$ is a closed path starting and ending at $B$ and containing no loops, then $\f F(l)=B^{\circ k},$ $k\geq 1,$ implying that 
 the image of $l$ under the homomorphism $\Phi_B$ from Theorem \ref{fg} is the unit element. 
Further, if $l$ contains a loop, then either 
$$\f F(l)=U_{kr}\circ \dots \circ U_{i+1}\circ \nu \circ  U_{i}\circ \dots \circ U_1,$$ 
 or
$$\f F(l)=U_{kr}\circ \dots \circ U_{i+1}\circ (\nu^{-1}\circ B_i) \circ  U_{i}\circ \dots \circ U_1$$
for some   $k\geq 1,$ $0\leq i \leq r-1,$ and $\nu \in Aut(B_i)$. Therefore, by Lemma \ref{le2} and Lemma \ref{sem0},   
either  
$$ \f F(l)=\mu \circ B^{\circ k},$$
or $$\f F(l)=\mu \circ B^{\circ (k+1)}$$
for some $\mu \in Aut(B).$  
 Finally, if 
$l$ contains several loops, then 
repeatedly using Lemma \ref{le2} and Lemma \ref{sem0} we conclude that 
$$\f F(l)=\mu \circ B^{\circ s}$$ for some $\mu \in Aut(B)$ and $s\geq 1.$
Thus, $G_B=Aut_G(B)$. \qed

\bc Let $B$  be a non-special  rational function of degree at least two such that 
$G_B$ is strictly larger than $Aut_G(B).$ Then there exists $A\sim B$ such that either $A$ can be represented 
as a composition of two commuting rational functions of degree at least two, or $A$ has more than one 
class of maximal decompositions. 
\ec
\pr By Theorem \ref{tde}, it is enough to show  that if any $A\sim B$ has a unique equivalence class of maximal decompositions and cannot 
be represented 
as a composition of two commuting rational functions of degree at least two,
then  for the function $B$ the both conditions defining generically decomposable rational functions are satisfied.
For the first condition this is obvious. For the second condition this is also true. 
Indeed, if
say $B_0=B$ is conjugate to $B_i$ and $\mu$ is a M\"obius transformation such that
$$ (U_r\circ \dots \circ U_{i+1})\circ   (U_{i}\circ\dots \circ U_1)= \mu\circ (U_{i}\circ\dots \circ U_1)\circ (U_r\circ \dots \circ U_{i+1})\circ \mu^{-1},$$ then   
for the functions 
$$N=\mu\circ (U_{i}\circ\dots \circ U_1), \ \ \  M= (U_r\circ \dots \circ U_{i+1})\circ \mu^{-1}$$
the equality 
\be \l{how} B=M\circ N=N\circ M \ee
holds. \qed

Notice that whenever $B$ {\it is} a composition of two commuting rational functions of degree at least two, the group $G_B$ is strictly larger than $Aut_G(B)$. Indeed, equality \eqref{how} implies easily that the functions $N$ and $M$ belong to $C_B$. Moreover, their images in  $G_B$ are not trivial and 
do not belong to $Aut_G(B)$, since $$1<\deg M< \deg B, \ \ \ 
1<\deg N< \deg B.$$ 
 In particular, if $B=T^{\circ s}$, where $s>1,$ the group $G_B$ contains a cyclic group of order $s$ whose intersection with $Aut_G(B)$ is trivial.

Finally, notice that the group $G_B$ can be strictly larger than $Aut_G(B)$ even if 
$B$ is not a composition of commuting functions, and 
 that the relation $A\sim B$ does not imply in general the equality $Aut_G(A)\cong Aut_G(B)$ 
(see Subsection 6.3).

\subsection{The group $G_B$ for polynomial $B$.}
Before stating the theorem describing groups $G_B$ for polynomial $B$ let us recall several results. 

First,
for a non-special polynomial $B$ of degree at least two, the set $C_B$ consists of {\it polynomials}. 
Indeed, \eqref{comm} yields that
\be \l{medva} B^{-1}(X^{-1}\{\infty\})=X^{-1}\{\infty\},\ee implying that 
$X^{-1}\{\infty\}$ contains at most two points. Furthermore, 
considering instead of $B$ and $X$ the  functions
\be \l{cha} X\rightarrow \mu\circ X\circ \mu^{-1}, \ \ \ \ B\rightarrow \mu\circ B\circ \mu^{-1}\ee  
for a convenient M\"obius transformation $\mu$, without loss of generality one can assume 
that either $X^{-1}\{\infty\}=\{\infty\}$, or $X^{-1}\{\infty\}=\{\infty,0\}.$ 
In the first case $X$ is a polynomial. On the other hand, in the second case \eqref{medva} implies that    
$B$  is conjugate to $z^n$,  in contradiction with the assumption that $B$ is not special. 

Secondly, the symmetry group $Aut(B)$ of a non-special polynomial $B$ of degree at least two is cyclic. Indeed, unless $B$  is conjugate to $z^n$,
for any $\mu\in Aut(B)$ necessarily  $\mu^{-1}\{\infty\}=\{\infty\}$, 
implying that $\mu$ is a polynomial. By a polynomial conjugation, we always can assume that the coefficient 
of $z^{\deg B-1}$ is zero, and it is clear that $\mu=az+b$ may commute with such $B$ only 
if $b=0$. Furthermore, it is easy to see that $Aut(B)$ is a cyclic 
rotation group of order $n$, 
where $n$ is the maximal number such that 
$$B=zR(z^n)$$ for some polynomial $R$.

 Thirdly, a polynomial $B$ is special if and only if $B$ is  conjugate to $z^{n}$ or $\pm T_n,$ since it is well known that a polynomial  cannot be a Latt\`es map.

In addition, we will need the following result (see \cite{pj}, Theorem 1.3).

\bt \l{enot}
Let $A$ and $B$ be fixed non-special polynomials of degree at least two, and  let $\f E(A,B)$ be the set of all polynomials of degree at least two $X$ such that $A\circ X=X\circ B$.
Then, either $\f E(A,B)$ is empty, or there exists $X_0\in \f E(A,B)$ such that a polynomial $X$ belongs to 
$\f E(A,B)$ if and only if $X= \widehat A\circ X_0$ for some polynomial  $\widehat A$ commuting with $A.$ \qed
\et

Recall that a group $G$ is called {\it metacyclic} if it has a normal cyclic subgroup $H$ such that $G/H$ is a cyclic group.

\bt \l{polo} Let $B$ be a polynomial of degree at least two not conjugate to $z^n$ or $\pm T_n$, $n\geq 2.$ Then the group $G_B$ is metacyclic.
\et 
\pr Applying Theorem \ref{enot} for $A=B$ and arguing as in Lemma \ref{ch1}, we see that any rational function $X$ that belongs to 
$C_B=\f E(B,B)$ has the form $X=\mu \circ X_0^{\circ l}$, where 
$\mu\in Aut(B)$ and $l\geq 1.$ In particular,  $ B=\mu \circ X_0^{l_0}$ for some $l_0\geq 1$ and $\mu \in Aut(B).$ 
 Moreover, the degree of any element of $C_B$ is a power of 
$d_0=\deg X_0$, and for $l\geq 0$ the subset of elements of degree $d_0^l$  coincides with the set 
$S_{1,l}=\{\mu \circ X_0^l\, \vert \, \mu \in Aut(B)\}.$  

Let us observe now that if \be \l{esz} X_0^{\circ l}\circ \mu_1=
X_0^{\circ l}\circ \mu_2,\ee where  $\mu_1, \mu_2 \in Aut(B),$ then $\mu_1=\mu_2.$ Indeed,  \eqref{esz}
implies that 
$$X_0^{\circ l}\circ (\mu_1\circ \mu_2^{-1})=X_0^{\circ l}.$$ Therefore, since 
$B^{\circ l}=\nu \circ X_0^{\circ (l_0l)}$ for some $\nu \in Aut(B),$
$$B^{\circ l}\circ (\mu_1\circ \mu_2^{-1})=B^{\circ l},$$ implying that $\mu_1=\mu_2.$ Thus, for $l\geq 0$ the set 
$S_{2,l}=\{ X_0^l\circ\mu  \, \vert \, \mu \in Aut(B)\}$ has the same cardinality as the set $S_{1,l}.$ Since $S_{2,l}$ is contained in $C_B$, this implies that $S_{1,l}=S_{2,l}.$   

The above analysis shows that 
the right cosets of $Aut_G(B)$ in $G$ have the form 
$$\b X_0^{l}Aut_G(B),   \ \ \ \ \ \ 0\leq l <l_0,$$
the left cosets  have the form $$Aut_G(B)\b X_0^{l},   \ \ \ \  \ \ 0\leq l <l_0,$$
and any right coset of $Aut_G(B)$ in $G$ is a left coset. 
Thus, $Aut_G(B)$ is a normal subgroup in $G_B$, and 
the group  $G_B/Aut_G(B)$ is a cyclic group of order $l_0$ generated by $\b X_0$. Since $Aut(B)$ is also  a cyclic group, we conclude that the group $G_B$ is metacyclic.
 \qed

Notice that Theorem \ref{polo} can be deduced from
the Ritt theorem (\cite{r}, \cite{r3}) saying that any commuting non-special polynomials $X$ and $B$ can be represented in the form 
\eqref{vot}. Nevertheless, the Ritt theorem does not imply  Theorem \ref{polo}  immediately, since 
$R$ in \eqref{vot} a priori depends on $X$, and the further analysis is needed.

\subsection{The group $G_B$ for the Ritt example.}
Let $B$ be a rational function of degree at least two. Denote by $\widehat{Aut}(B)$ the group consisting of M\"obius transformations $\mu$ such that 
$$B\circ \mu=\nu\circ B$$ for some M\"obius transformations $\nu$. 
Like the group $Aut(B)$, the group  $\widehat{Aut}(B)$ is a finite rotation group
of the sphere (see \cite{fin}, Section 4). More generally, denote by $\widehat C_B$ the set of rational functions $X$ such that 
$$B\circ X=Y\circ B$$ for some rational function $Y.$ 
Clearly, $Aut(B)$ is a subgroup  $\widehat{Aut}(B)$, and  $C_B\subseteq \widehat C_B.$

Let 
$$V=\frac{z^2+2}{z+1},\ \ \ U=\frac{z^2-4}{z-1},\ \ \ \mu=\varepsilon z,$$ where $\varepsilon^3=1$. In the paper \cite{r}, Ritt showed that 
the rational functions 
$$B=V\circ U, \ \ \ \ X=V\circ \mu \circ U$$ commute but no one of them is a rational function of the other. In particular, this implies that there is no $R$ such that 
\be \l{tutb} B=\mu_1\circ R^{\circ l_1}, \ \ \ X=\mu_2\circ R^{\circ l_2}\ee
for some M\"obius transformations $\mu_1,$ $\mu_2$ and $l_1,l_2\geq 1.$ More generally, for 
any function $C$ such that $C(\varepsilon z)=\varepsilon C(z)$
 the functions 
$$B'=V\circ C\circ U, \ \ \ \ X'=V\circ \mu \circ C \circ U$$
commute but no one of them is a rational function of the other.

The Ritt statement follows from the following more general observation.

\bl \l{rl+} Let $W \in C_{U\circ V}$ but $W \not\in \widehat{C}_V$. Then 
the functions $V\circ U$ and $V\circ W\circ U$ 
commute but the latter is not a rational function of the former. Furthermore, the same 
conclusion holds for 
the functions $V\circ C\circ U$ and $V\circ W \circ C \circ U,$ where $C$ is any function commuting with $W$.
\el
\pr Indeed, we have:
 $$(V\circ C\circ U)\circ (V\circ W\circ C\circ U)=V\circ C\circ (U\circ V\circ W)\circ C\circ U=$$ $$
=V\circ C\circ ( W\circ U\circ V)\circ C\circ U=(V\circ C\circ  W\circ U)\circ (V\circ C\circ U)=$$
$$=(V\circ W\circ  C\circ U)\circ (V\circ C\circ U).$$
On the other hand, if 
\be \l{im} V\circ W\circ C\circ U=R\circ V\circ C\circ U\ee  for some rational function $R$, then 
$$V\circ W=R\circ V,$$ in contradiction with the assumption that 
 $W \not\in \widehat{C}_V$.  \qed

The Ritt statement is obtained from Lemma \ref{rl+} for $W=\mu$.  Indeed,   
$$U\circ V={\frac {z \left( {z}^{3}-8 \right) }{ \left( {z}^{3}+1 \right) 
  }},$$ implying that $\mu \in Aut(U\circ V)$. On the other hand, 
the assumption that 
\be \l{gop} V\circ \mu=\nu\circ V\ee for some 
M\"obius transformation $\nu$ leads to a contradiction. Namely, 
\eqref{gop} implies that $\nu(\infty)=\infty$. Therefore,  $\nu=az+b,$ $a,b\in \C,$ and hence  
if \eqref{gop} holds, then the functions $V$ and 
$$V\circ \mu=\frac{\v^2z^2+2}{\v z+1}$$ have the same set of poles.
However, this is not true.

Let us calculate the group $G_B$. Using again a computer assistance 
one can check that the function 
$$B=V\circ U={\frac {{z}^{4}-6\,{z}^{2}-4\,z+18}{ \left( {z}^{2}+z-5 \right) 
 \left( z-1 \right) }}
$$ has four critical values and the corresponding permutations in $Mon(B)$ can be identified with  
the permutations $(1 3),$ $(12)(34)$, $(1 3),$ and $(12)(34)$ in $S_4$, while the function 
$$B_1=U\circ V={\frac {z \left( {z}^{3}-8 \right) }{ \left( {z}^{3}+1 \right) 
  }}$$ has three critical values and the corresponding permutations in $Mon(B_1)$  can be identified with  
$(12)(34),$ $(13)(24),$ and $(14)(23).$
In particular, $B_1$ and $B$ are not conjugate  since they have a different number of critical 
values. Moreover, one can check that the group $Aut(B)$ is trivial while $Aut(B_1)$ is a cyclic group of order three 
generated by $\mu.$

\begin{figure}[htbp]
\begin{center}
\epsfig{file=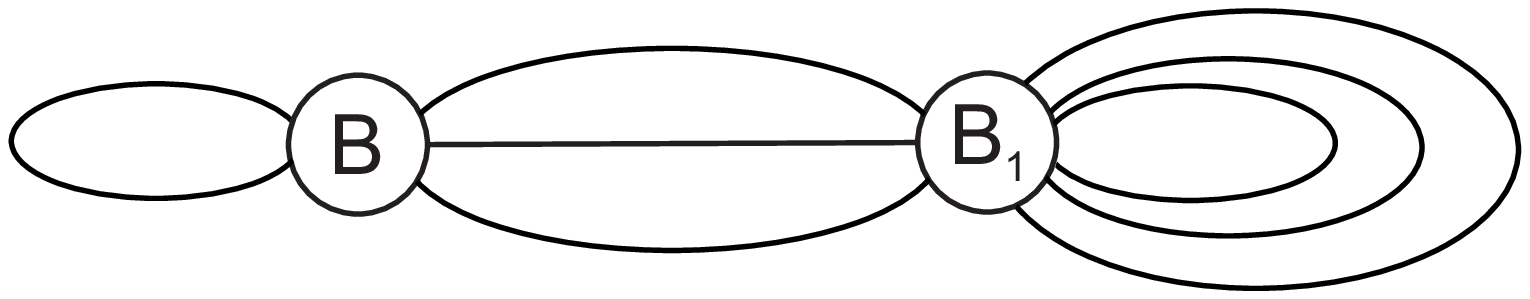,width=7.5cm}
\caption{ }
\label{fig4}
\end{center}
\end{figure}

It is easy to see that  $Mon(B)$ has a unique imprimitivity system
$\{1,3\},\{2,4\}$, corresponding to the decomposition $B=V\circ U$ while 
$Mon(B_1)$ has three imprimitivity systems 
 $$\{1,3\},\{2,4\}, \ \ \ \ \{1,2\},\{3,4\}, \ \ \ \ \ \{1,4\},\{2,3\}.$$ 
 corresponding to the decompositions
$$B_1=U\circ V, \ \ \ \ \ \ \ \ \  B_1=(\mu^{-1}\circ U)\circ (V\circ \mu), \ \ \ \ \ \ \ \ \  B_1=(\mu^{-2}\circ U)\circ (V\circ \mu^2).$$
Summing up, we see that the graph $\Gamma_B$ has the form shown on 
 Fig. \ref{fig4}, 
where the edges connecting $B$ and $B_1$  correspond to the solutions
$$B=( V\circ\mu^{i-1})\circ (\mu^{-(i-1)}\circ U), \ \ \ \ \ 
 B_1=(\mu^{-(i-1)}\circ U)\circ (V\circ \mu^{i-1}), \ \ \ \ \ 1\leq i \leq 3,$$
of system \eqref{sys1}, the loops attached to $B_1$  correspond to the solutions
$$B_1=(\mu^{-(i-1)}\circ B_1)\circ \mu^{i-1}=\mu^{i-1}\circ (\mu^{-(i-1)}\circ B_1),
\ \ \ \ \ 1\leq i \leq 3,$$
and the loop attached to $B$ corresponds to the solution \eqref{ewq0}.

\begin{figure}[htbp]
\begin{center}
\epsfig{file=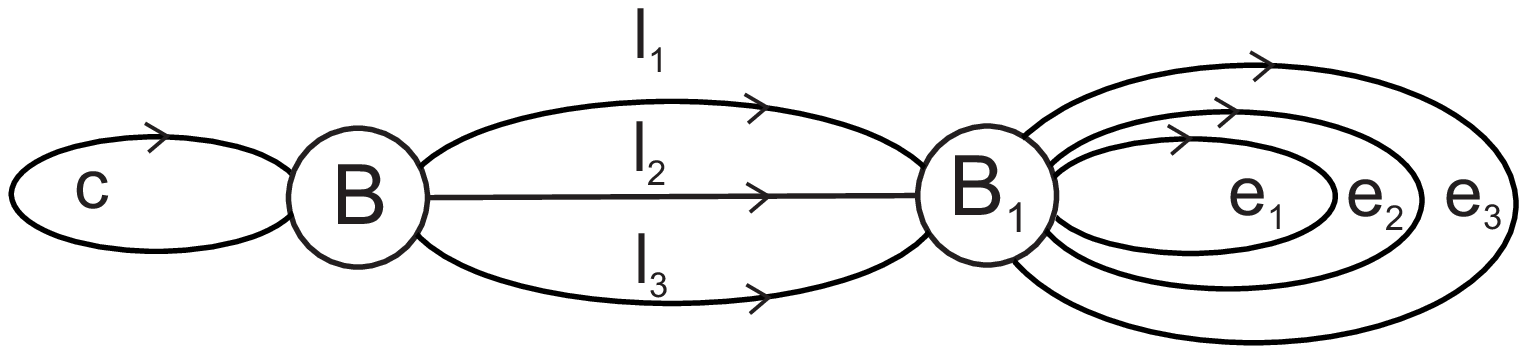,width=7.5cm}
\caption{ }
\label{fig6}
\end{center}
\end{figure}

The fundamental group of $\Gamma_B$ can be easily calculated by the well known method 
using the spanning tree (see e. g. \cite{still}, Section 4.1.2). 
Namely, choosing a fixed orientation on each of edges of $\Gamma_B$ as it is shown 
on Fig. \ref{fig6},  
and considering the edge $l_1$ together with vertices $B$ and $B_1$ as the spanning tree, 
we see that $\pi_1(\Gamma_B,B)$ is a free group of rank 6 generated by the paths
$$c, \ \ \ \ l_1^{-1}l_i, \ \ \ 2\leq i \leq 3, \ \ \ l_1^{-1}e_jl_1, \ \ \ 1\leq j\leq 3,$$
implying that the group $G_B$ is generated by the images of these paths under the 
map $\Phi_B$. Assuming that 
$$\f F(c)=z, \ \ \ \ \f F(e_i)=\mu^{i-1}, \ \ \ \ 1\leq i \leq 3,$$ 
we obtain
$$\f F(l_1^{-1}l_i)=V\circ\mu^{-(i-1)}\circ U, \ \ \ 2\leq i \leq 3, \ \ \  \f F(l_1^{-1}e_jl_1)=V\circ\mu^{j-1}\circ U, \ \ \ 1\leq j\leq 3,$$ implying that 
the images of the functions 
\be \l{thef} g_0=z, \ \ \ \ \ g_1=V\circ\mu \circ U, \ \ \ \ \ g_2=V\circ\mu^2 \circ U\ee under the map $\Phi_B$  generate the 
group $G_B.$ Since 
\be \l{edg} \deg g_1=\deg g_2=\deg B,\ee and
$$g_1\neq B, \ \ \ \ g_2\neq B, \ \ \ \ g_1\neq g_2,$$ it follows from Lemma \ref{eg} that 
$g_1,g_2,g_3$ represent different classes in $C_B/\underset{B}{\sim}$\,, so that 
$G_B$ has at least three elements. On the other hand, we have: 
$$g_1^{\circ 2}=g_2\circ B, \ \ \ \ \ g_2^{\circ 2}=g_1\circ B, \ \ \ \ \  g_1^{\circ 3}=g_2^{\circ 3}=B^{\circ 3},\ \ \ \ \ g_1\circ g_2=g_2\circ g_1=B^{\circ 2}.$$
Therefore, $G_B=\Z/3\Z$.

In turn, the set $C_B$ can be described as follows:  $X\in C_B$ if and only if either
$$X=B^{\circ j}, \ \ \ \ \ \ j\geq 0,$$ or
$$X=V\circ\mu \circ U \circ B^{\circ j},\ \ \ \ \ \ j\geq 0, $$ 
or 
$$X=V\circ\mu^2 \circ U \circ B^{\circ j}, \ \ \ \ \ j\geq 0.$$ 
Indeed, by Lemma \ref{eg}, it is enough to check that the functions \eqref{thef} 
are not rational functions in $B$. 
Assume say that $ g_1=R\circ B.$ Then it follows from \eqref{edg} 
that 
$R$ is a M\"obius transformation. Moreover,  $R\in Aut(B)$ by Lemma \ref{per}. However, since 
$Aut(B)$ is trivial and $g_1\neq B,$ this is impossible. 

Notice that since $G_B\cong G_{B_1}$ by Theorem \ref{var} and $Aut_G(B_1)=\Z/3\Z,$ we have: $$G_{B_1}=Aut_G(B_1)=\Z/3\Z.$$
Notice also that since $G_B\cong G_{B_1}$, the non-triviality of  $Aut(B_1)$ already implies the non-triviality of  $G_B$. Moreover, since $B$ has no automorphisms, we can conclude that the set $C_B$ contains functions of degree greater than one that are not iterates of $B.$

\subsection{The group $G_B$ for $B=-2z^2/(z^4+1)$.}
\begin{figure}[htbp]
\begin{center}
\epsfig{file=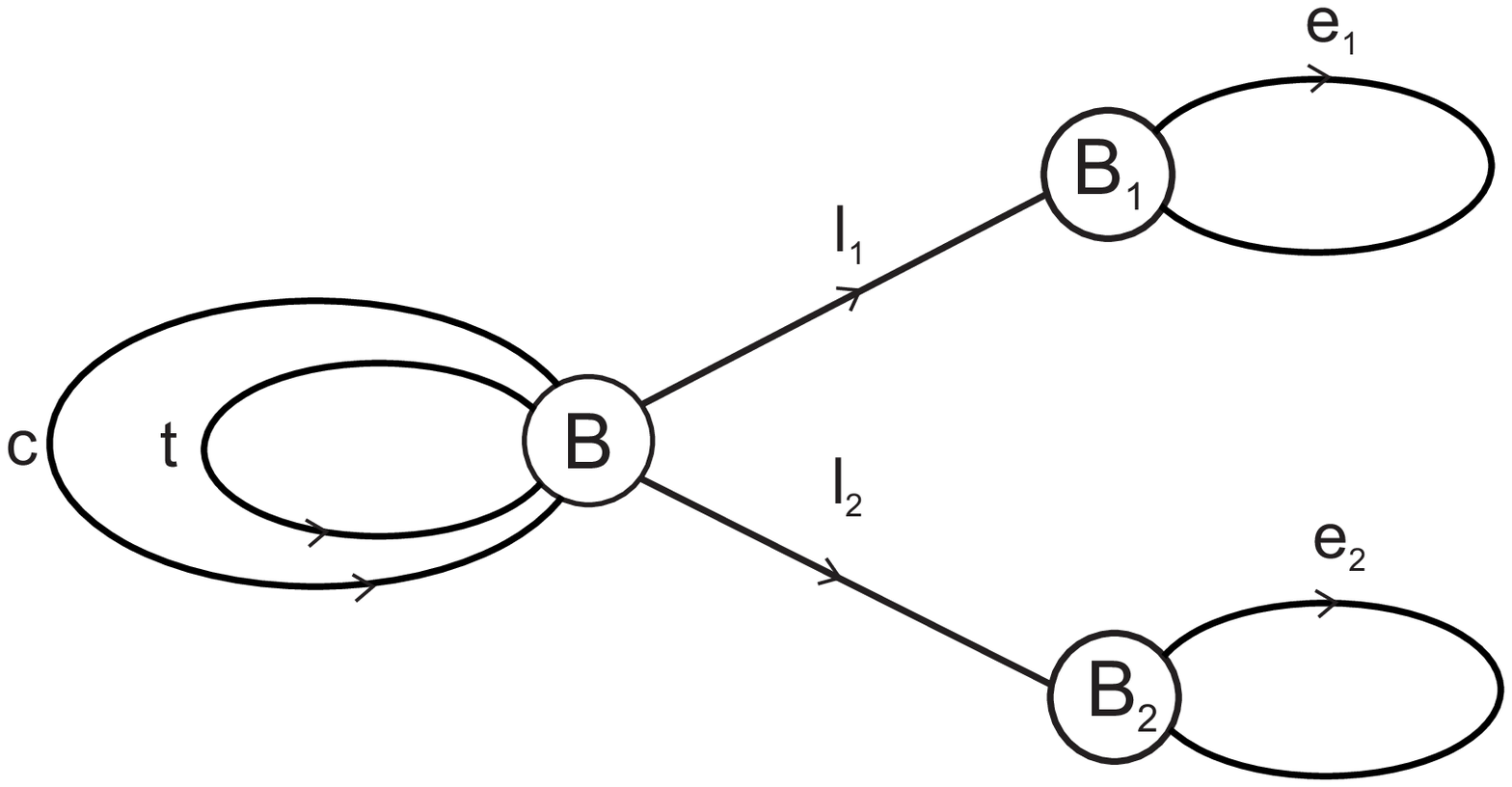,width=5.9cm}
\caption{ }
\label{fig7}
\end{center}
\end{figure}

Since equality \eqref{raka} implies that 
the function \be \l{db} W=\frac{z^2-1}{z^2+1}\ee commutes with $B$, the group $G_B$ 
clearly contains a cyclic group of order two generated by $\b W$. Moreover, 
it is easy to see that in fact $G_B=\Z/2\Z$.
Indeed, 
providing edges of the graph $\Gamma_B$ with   orientations shown   on Fig. \ref{fig7},  
we see that $\pi_1(\Gamma_B,B)$  
is a free group of rank 4 with  generators
$$c, \ \ \ t, \ \ \ l_i^{-1}e_il_i, \ \ \ i=1,2,$$ and 
assuming that 
$$\f F(c)=\f F(e_1)=\f F(e_2)=z,\ \ \ \ \ \f F(t)=W,$$ we see that 
$G_B$ is generated by the $\b W.$ Similarly, one can conclude that 
$G_{B_1}$ is genera\-ted by  $\b X,$
where \be \l{te} X=\f F(l_1tl_1^{-1})=\frac{z^2+1}{z}\circ \frac{z^2-1}{z^2+1}\circ -\frac{2}{z^2-2}.\ee

The above functions $B_1$ and $X$ provide an example of commuting rational functions similar to the one 
constructed by Ritt.
Namely, set 
 $$V=\frac{z^2+1}{z}, \ \ \ U=-\frac{2}{z^2-2}.$$ Then $W$ commutes with $U\circ V=W^{\circ 2},$ but 
 $W \not\in \widehat{C}_V$. Indeed, assume the inverse, and let $S$ be the rational function defined by any of the sides  of the equality
\be \l{gopa} \frac{z^2+1}{z}\circ \frac{z^2-1}{z^2+1}=R\circ \frac{z^2+1}{z}, \ee
where $R\in \C(z)$.  Then substituting $z$ by $\frac{1}{z}$ in the right side 
of \eqref{gopa} we obtain that $S\circ \frac{1}{z}=S.$ However, substituting $z$ by $\frac{1}{z}$ in the left side we obtain 
$$S\circ \frac{1}{z}=\frac{z^2+1}{z}\circ -\frac{z^2-1}{z^2+1}=-S.$$ The contradiction obtained shows that $W \not\in \widehat{C}_V$. Therefore,  
by Lemma \ref{rl+}, the rational function 
$$X=V\circ W\circ U$$ commutes with $B_1=V\circ U$, but is not a rational function in $B_1$.
Notice that
in distinction with the Ritt example 
the  non-triviality of $G_{B_1}$ is explained by the existence 
in the class $[B_1]$  of a function that is an iterate.

\vskip 0.2cm
\noindent{\bf Acknowledgments}.
The author is grateful to  
the Max-Planck-Institut fuer Mathematik for the hospitality and the support.

\bibliographystyle{amsplain}

\end{document}